\documentclass[onesided,ref,12pt]{amsart}
\usepackage{amscd,amsfonts,amssymb,amstext,amsthm,eufrak,yfonts}
\usepackage{mathbbol,mathrsfs,stmaryrd}
\usepackage[cmtip,arrow]{xy}
\usepackage{pb-diagram,pb-xy} 
\usepackage{graphics} 
\newcommand{\comment}[1]{}
\newcommand{\pfend}{\hfill{$\boxempty$}}

\newtheorem{definition}{\sc Definition}[section]
\newtheorem{theorem}{\sc Theorem}[section]
\newtheorem{proposition}{\sc Proposition}[section]

\title{General Concepts of Graphs} 
\author{}
\date{} 

\begin{document} 
\maketitle 
\vskip 5mm 
\begin{center} 
\textbf{Dedicated to the 65th anniversary of 
Professor Ulrich Knauer}\\[5mm] 
Sheng Bau\\[3mm] 
Center for Discrete Mathematics,\\  
Fuzhou University, Fuzhou, China\\ 
\texttt{dimacs1@fzu.edu.cn}\\[2mm]   
and\\[2mm] 
Schoold of Mathematical Sciences,\\ 
Government College University, Lahore, Pakistan 
\end{center}  
\vskip 3mm 
\hrule 
\vskip 3mm 

\begin{abstract} 
A little general abstract combinatorial nonsense 
delivered in this note is a presentation of some 
old and basic concepts, central to discrete mathematics, 
in terms of new words. The treatment is from 
a structural and systematic point of view. 
This note consists essentially of definitions and summaries. 
\end{abstract} 
\vskip 3mm 
\hrule 
\tableofcontents 
\hrule 
\newpage 
\section{Generalization and Specialization}  
\label{sec:genspec} 
There had been a general trend of generalization 
throughout 19th and the first half of the 20th centuries 
in mathematics. Then in the second half of the 20th 
century, a prevailing trend of specialization occurred. 
The latter trend is characteristic of almost exclusive 
emphasis on the immediate solutions of specific problems, 
especially if they were openly proposed by another 
mathematician, usually a famous one. This was partly 
encouraged by the ad hoc aim of seeking \textsl{academic excellence} 
by journals and by their authors. If a paper is on 
some systematic ground work which entails in proper 
generalization and exposure of important relations 
between fundamental mathematical concepts, 
it might easily be mistaken as ``no novelty'' and 
will be appropriated an instant rejection, sometimes 
even without a proper refereeing procedure. 

This note consists only of defintions of basic concepts. 
With these, I wish to emphasize that mathematics is a unity. 
Individual activities and areas of activity are related in 
a vitally organic manner. 
A monopoly of dissection of the body of mathematics 
is clearly not always beneficial to the health and life of 
mathematics. 

In support of the view that mathematics is 
an organic unity and many branches of our 
science are vitally related, I would like 
to point out the intensive and extensive interactions 
between algebra and combinatorics 
as in \cite{godsil:1993,godsilr:2001} 
and the literature therein, 
between probability and combinatorics 
(\cite{alons:1992,bollobas:1985} and references therein), 
between topology and combinatorics \cite{grosst:2001,mohart:2001}, 
and more recently, between analysis and combinatorics 
\cite{taov:2006}. Mathematics is, after all, \textsl{not} 
an exclusive instrument of some small number of 
executives in the community for ad hoc 
\textsl{academic excellence}. 

The concepts of graphs and some of their generalizations 
are included in Section \ref{sec:relsys}. These will be 
specialized to algebraic objects in Section \ref{sec:mathobj}. 
Different types of categories of graphs 
will be reviewed in Section \ref{sec:catgraphs}, 
while the concept of graph invariants will 
be clarified in Section \ref{sec:invar}. 
Inductive sets arise naturally as graphs 
in Section \ref{sec:recurs}. 
Transformation graphs arise naturallly in mathematics. 
Generality of the concept of transformation graphs will 
be considered in Section \ref{sec:transfg}. 

\section{Graphs and Relational Systems} 
\label{sec:relsys} 
\begin{definition}{\rm 
Let $V$ and $E$ be sets with $V\cap E =\emptyset$. 
Then a mapping $G : E\rightarrow V\times V$ 
is called a \textit{graph}. This may be given by 
\[
\begin{diagram} 
\node{E}\arrow{e,t}{G}\node{V\times V} 
\end{diagram}
\] 
} 
\label{def:graph} 
\end{definition} 
This definition captures the essence of the concept 
of a graph precisely and with proper generality, as it includes 
the concepts of both finite and infinite, both simple 
and nonsimple, both directed and undirected graphs. 

A special case, when $G$ is an injective mapping, 
of Definition \ref{def:graph} may be presented 
in terms of a binary relation. 

\begin{definition}{\rm 
Let $V$ be any set and $E\subseteq V\times V$ 
be a binary relation over $V$. 
Then the ordered tuple $G = (V, E)$ is called a \textit{graph}. 
} 
\label{def:graphasbinrel} 
\end{definition} 

A concept of this generality is native to mathematics. 
A directed graph with multiple edges that 
are weighted may be presented by 
\[
\begin{diagram} 
\node{W}\node{E}\arrow{w,t}{w}\arrow{e,t}{G}\node{V\times V} 
\end{diagram} 
\] 
or 
\[
\begin{diagram} 
\node{}\node{W}\arrow{s,b}{w}\\ 
\node{\phantom{V()}(V,}\node{E)}
\end{diagram}
\] 
where $W\subseteq\mathbb{R}$ and 
$w$ is a nonnegative function. 

If $V$ is finite then $G$ is called \textit{finite}; 
if $E$ is irreflexive, that is, 
if for each $a\in V$, $(a, a)\not\in E$, 
then $G$ is \textit{loopless}; 
if $E$ is symmetric then the tuples in $E$ may be considered 
subsets of cardinality $2$ and hence $G$ is \textit{undirected}; 
if $E$ is antisymmetric then $G$ is \textit{oriented}; 
if $E$ is antisymmetric and transitive then $G$ is a 
\textit{partial order}, an immediate specialization. 

\begin{proposition} 
Every partially ordered set is an oriented graph. 
\label{prop:poset} 
\end{proposition} 
Hence, every lattice is an oriented graph. 

Denote 
\[D(V) =\{(a, a) : a\in V\}\] 
and call $D(V)$ the \textit{diagonal} of 
$\times$, as usual. Note that the diagonal 
is also a binary relation on $V$. 
Let $E$ be a symmetric binary relation on $V$. 
The \textit{symmetric closure} $s(E)$ of 
$E$ is 
\[s(R) := E\cup \{(b, a) : (a, b)\in E\}.\] 
Then the natural projection 
\[p_s : s\left[V\times V - D(V)\right]\rightarrow [V]^2 :=
\left(\begin{array}{c}V\\ 2\end{array}\right)\]
of the symmetric closure 
forgets the \textsl{order} of ordered tuples 
and maps ordered tuples to subsets of cardinality $2$. 
If 
\[G : E\rightarrow [V]^2\] 
is a bijection, then $G$ is called the 
\textit{complete} graph of order $|V|$ and 
is denoted $G = K_{|V|}$. Any simple graph 
$G$ with $V(G) = V$ is a \textit{spanning} subgraph 
of $K_{|V|}$ since 
\[ 
\begin{diagram} 
\node{}\node{E(K_{|V|})}\arrow{s,t}{K_{|V|}}\\
\node{E(G)}\arrow{ne,t}{K^{-1}_{|V|}\cdot G}
\arrow{e,t}{G}\node{[V]^2} 
\end{diagram} 
\] 
Since $K_{|V|}$ is a bijection, and $G$ is injective, 
hence 
\[K^{-1}_{|V|}\cdot G : E(G)\rightarrow E(K_{|V|})\] 
is an injective mapping (i.e., an embedding). 
That is, the spanning subgraph 
\[G = K_{|V|}\cdot K^{-1}_{|V|}\cdot G.\] 

For a graph $G$, the morphism \textit{underlying undirected} 
graph provides a forgetful functor from the category of graphs 
to the category of undirected graphs, presented 
by the following diagram. The forgetful 
functor $U(G)$ forgets the directions of edges. 
If a graph is denoted $\overrightarrow{G}$ then 
its underlying undirected graph $U(\overrightarrow{G})$ 
may be denoted by $G$. 
\[ 
\begin{diagram} 
\node{E}\arrow{e,t}{\overrightarrow{G}}\arrow{se,b}{G}
\node{V\times V}\arrow{s,b}{p'_s}\\ 
\node{}\node{D(V)\cup [V]^2}
\end{diagram} 
\] 
where $p'_s$ is an extension of $p_s$ over to 
$D(V)\cup [V]^2$. 
A \textit{simple graph} may be presented 
by 
\[ 
\begin{diagram} 
\node{E}\arrow{e,t}{G}\node{[V]^2.}
\end{diagram} 
\] 

The following diagram presents mappings 
of edges of a graph to their heads ($p_2G$) 
and tails ($p_1G$). 
\[ 
\begin{diagram} 
\node{}\node{E}\arrow{sw,t}{p_1G}
\arrow{s,t}{G}\arrow{se,t}{p_2G}\\ 
\node{V}\node{V\times V}
\arrow{w,b}{p_1}\arrow{e,b}{p_2}
\node{V} 
\end{diagram} 
\] 
This gives information about incidence. 
The \textit{incidence matrix} $M_G$ 
of a loopless graph $G$ is given by 
\[M_G : V\times E\rightarrow \{-1, 0, 1\}.\] 
Thus the incidence matrix is a matrix 
whose rows are indexed by $V$, columns by $E$ 
and entries in $\{-1, 0, 1\}$ which reveals 
the manner of incidence of $v\in V$ and $e\in E$ 
as given below. 
\[M_G(v, e) =\left\{
\begin{array}{ll} 
-1, & G(e) = (v, w),\; w\in V\\ 
\phantom{-}0, & v\not\in e\\ 
\phantom{-}1, & G(e) = (u, v),\; u\in V 
\end{array}
\right. 
\] 

Let $X\subseteq E$ where $\iota : X\rightarrow E$ 
is the inclusion. Then the subgraph 
$G|_X\subseteq G$ \textit{induced} by $X$ 
is presented by the diagram 
\[ 
\begin{diagram} 
\node{X}\arrow{e,t}{\iota}\arrow{se,b}{G|_X}
\node{E}\arrow{s,b}{G}\\ 
\node{}\node{V\times V} 
\end{diagram} 
\] 
That is, $G|_X = G\cdot\iota$. 
For $S\subseteq V$ with inclusion 
$\eta : S\rightarrow V$, the subgraph 
$G|_S$ \textit{induced} by $S$ is presented as 
\[ 
\begin{diagram} 
\node{E}\arrow{e,t}{G}\node{V\times V}\\ 
\node{G^{-1}(S\times S)}\arrow{n,t}{\eta'}
\arrow{e,b}{G|_S}\node{S\times S}\arrow{n,b}{\eta_{S\times S}} 
\end{diagram} 
\] 
where $\eta_{S\times S} : S\times S\rightarrow V\times V$ 
is the inclusion mapping induced by $\eta$, and 
$\eta' : G^{-1}(S\times S)\rightarrow E$ 
is the natural inclusion. 
The diagram commutes: $\eta_{S\times S}\cdot G|_S = G\cdot\eta'$. 

In this note, some further generalization 
will be considered before moving to the points of 
specializations. There are obviously two 
directions in which the concept of a graph 
as given in Definition \ref{def:graph} may 
be generalized, generalization at the head 
or at the tail of the arrow in the diagram. 

\begin{definition}{\rm 
Let $V$ be a set and let $E = (E_1, \cdots , E_s)$ 
be a collection of sets. If 
\[G = (G_1, \cdots , G_s)\] 
and 
\[G_i : E_i\rightarrow V\times V,\; i = 1, \cdots , s\] 
then 
\[G : (E_1, \cdots , E_s)\rightarrow V\times V\] 
is called a \textit{graph system}. 
A graph system may be presented by the diagram 
\[ 
\begin{diagram} 
\node{E_1}\arrow{se,t}{G_1}\\ 
\node{\vdots}\node{V\times V}\\ 
\node{E_s}\arrow{ne,b}{G_s}
\end{diagram} 
\] 
} 
\label{def:graphsys} 
\end{definition} 

This seems to be an obvious generalization of the concept 
of a graph at its generality as given in Definition \ref{def:graph}, 
though, as may be seen directly from the diagram, 
it may be fully realized as a collection 
of graphs over the same set $V(G) = V$. 

The second generalization is at the head 
of the arrow in the diagram of Definition 
\ref{def:graph}. 

Denote 
\[V^m :=\underbrace{V\times\cdots\times V}_{m\;\textrm{fold}}.\] 

\begin{definition}{\rm 
Let $V$ and $E$ be sets. 
\[G : E\rightarrow\underbrace{V\times\cdots\times V}_{m\;\textrm{fold}}\] 
is called a \textit{set system} or a \textit{hypergraph}. 
A set system (hypergraph) may be presented by the diagram 
\[ 
\begin{diagram} 
\node{E}\arrow{e,t}{G}\node{V^m}
\end{diagram} 
\] 
} 
\label{def:setsys} 
\end{definition} 

This is a proper generalization. Note that this 
is more general than the usual a set system or 
a hypergraph since at the head of the arrow 
is a cartesian product, instead of 
\[\bigcup\limits^m_{k=1}[V]^m 
:=\{S\subseteq V : |S|\leq m\}.\]  
For $R\subseteq V^m$, 
the \textit{symmetric closure} or 
$S_m$-\textit{closure} $s(R)$ is 
the quotient defined by the binary relation 
\[(v_{p(1)}, v_{p(2)}, \cdots, v_{p(m)})
\sim (v_1, v_2, \cdots , v_m)\] 
for all $p\in S_m$ the symmetric group 
on $\{1, \cdots , m\}$. 
That is, $s(R)$ is the quotient of the transitive action of 
the symmetric group $S_m$, then 
the usual set system (hypergraph) 
may be presented by the diagram 

\[ 
\begin{diagram} 
\node{E}\arrow{e,t}{G}
\node{\bigcup\limits^m_{k=1}\left(\begin{array}{c}V\\ k\end{array}\right)}
\end{diagram} 
\] 

The concept of a \textit{relational system} 
may be obtained by a generalization on 
both head and tail of the arrow in 
Definition \ref{def:graph}, and insisting 
that the mappings concerned are injections. 

\begin{definition}{\rm 
Let $V$ be a set and let $E = (E_1, \cdots , E_s)$ 
be a collection of sets. If 
\[G = (G_1, \cdots , G_s)\] 
and 
\[G_i : E_i\rightarrow\bigcup\limits^m_{k=1}[V]^m, \; i = 1, \cdots , s\] 
then 
\[G : (E_1, \cdots , E_s)\rightarrow\bigcup\limits^m_{k=1}[V]^m\] 
is called a \textit{relational system}. 
A relational system may be presented by the diagram 
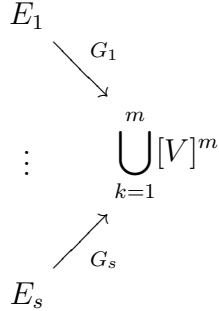
\[ 
\begin{diagram} 
\node{E_1}\arrow{se,t}{G_1}\\ 
\node{\vdots}\node{\bigcup\limits^m_{k=1}[V]^m}\\ 
\node{E_s}\arrow{ne,b}{G_s}
\end{diagram} 
\] 
} 
\label{def:relsys} 
\end{definition} 

If each $G_i$ in Definition \ref{def:relsys} 
is injective, then we have 

\begin{definition}{\rm 
Let $V$ be any set and let $E$ be a collection of 
relations over $V$. Then $G = (V, E)$ is called 
a \textit{simple relational system}. Let $k_1 < k_2 <\cdots < k_m$ 
and let the number of distinct $k_i$-ary relations in $E$ 
be $r_i$. Then the symbol $(k^{r_1}_1, k^{r_2}_2, \cdots , k^{r_m}_m)$ 
is called the \textit{type} of $G$. The integer $k_m$ 
is called the \textit{arity} of $E$ and hence of $G$.  
} 
\label{def:srelsys} 
\end{definition} 

Trivially, every hypergraph is a relational system. 

Now a proper specialization of the concept of a relational system. 
\begin{definition}{\rm 
A simple relational system $G$ with arity $3$ is called 
a \textit{ternary relational system}. 
} 
\label{def:ternrelsys} 
\end{definition} 

This is the definition that encompasses almost all 
mathematical objects. The first and the most important concept 
is the concept of a group. Since each binary operation is a ternary 
relation, \textsl{every group is a ternary relational system}. 
At the appropriate level of binary operations, algebraic objects 
(rings, principal ideal domains, division rings, and fields) 
have been studied extensively. Indeed, algebra represents 
one of the great successes in modern mathematics. 
At this point, there exists a rich possibility of 
specialization of the concept of relational systems, 
which certainly leads to an abundance of problems and questions 
including the investigation of graphs, groups 
and partially ordered sets. 

The general concept of a relational system 
is a very recent one \cite{helln:2004}, 
and very little is known about relational systems, 
while obviously there is a rich and extensive theory 
in the case where these relations are operations. 
As a suitable generality is now at hand, this is the point 
where a few specializations will be considered more formally. 
But, before this will be dealt with in Section \ref{sec:mathobj}, 
consider a few typical categories of graphs. 

\section{Categories of Graphs} 
\label{sec:catgraphs} 
There are many different categories of graphs 
(See \cite{kilpk:2001,kilpkm:2000,knauer:2006}). 
For a formal and comprehensive treatment of categories and 
functors, see \cite{maclane:1971}. In this section, 
some fundamental categories of graphs will be 
presented. This section is based on an excerpt 
from \cite{kilpk:2001,knauer:2006}. 

\begin{definition}{\rm 
\begin{enumerate} 
\itemsep 0pt\parskip 0pt 
\item 
The category of all graphs with 
graph homomorphisms as morphisms. 
This category is denoted by $\textswab{Gra}$. 
A \textit{homomorphism} $f : G\rightarrow H$ is 
a mapping with $xy\in E(G)\Rightarrow f(x)f(y)\in E(H)$. 
\item 
The category of all graphs with egamorphisms as morphisms. 
This category is denoted by $\textswab{EGra}$. 
An \textit{egamorphism} is a mapping $f : G\rightarrow H$ such 
that $xy\in E(G)\Rightarrow f(x)f(y)\in E(H)\cup V(G)$. 
\item 
The category of all graphs with comorphisms as morphisms. 
This category is denoted by $\textswab{CGra}$. 
A \textit{comorphism} is a mapping $f : G\rightarrow H$ such that 
$f(x)f(y)\in E(H)\Rightarrow xy\in E(G)$.  
\end{enumerate} 
}
\label{def:catgraphs} 
\end{definition} 

In each of the three categories, compositions 
and morphisms, respectively, 
obviously satisfy the categorical axioms for 
compositions. Note that one of many ways of 
definining the important concept of a \textsl{contraction} 
is that it is a preconnected egamorphism, meaning 
that it is an egamorphism (i.e., morphism of 
the category $\textswab{EGra}$ for which 
the preimage of each vertex induces a connected 
subgraph). This will be addressed by the author 
in another paper. 

The binary operations of graphs typically include 
products. Natural products in respective 
graph categories will now be reviewed. 

Consider the category $\textswab{S}$ of sets 
and mappings. Let $G_1, G_2\in\textswab{S}$. 
A pair $(G, (p_1, p_2))$ with 
$p_1 : G\rightarrow G_1$, $p_2 : G\rightarrow G_2$ 
is called (the categorical) \textit{product} of 
$G_1, G_2$ in $\textswab{S}$ if 
(1) $p_1, p_2$ are morphisms in $\textswab{S}$; and 
(2) $(G, (p_1, p_2))$ solves the universal problem: 
for all sets $H$ and for all mappings 
$f_1 : H\rightarrow G_1$, $f_2 : H\rightarrow G_2$ 
there exists a unique mapping $f : H\rightarrow G$ 
such that the diagram 
\[
\begin{diagram} 
\node{}\node{H}\arrow{sw,t}{f_1}\arrow{s,l}{f}\arrow{se,t}{f_2}\\ 
\node{G_1}\node{G}\arrow{w,t}{p_1}\arrow{e,t}{p_2}\node{G_2} 
\end{diagram}
\] 
commutes. 

\begin{theorem} 
$(G_1\times G_2, (p_1, p_2))$ is the product of 
$G_1$ and $G_2$ in $\textswab{S}$. 
\label{the:product} 
\end{theorem} 

A pair $((u_1, u_2), G)$ is 
called the \textit{coproduct} of $G_1, G_2$ in 
$\textswab{S}$ if $u_1 : G_1\rightarrow G$, 
$u_2 : G_2\rightarrow G$ are mappings such that 
for all sets $H$ and for all mappings 
$f_1 : G_1\rightarrow H$, $f_2 : G_2\rightarrow H$ 
there exists exactly one mapping 
$f : G\rightarrow H$ such that the diagram 
\[
\begin{diagram} 
\node{G_1}\arrow{e,t}{u_1}\arrow{se,b}{f_1}
\node{G}\arrow{s,l}{f}\node{G_2}\arrow{w,t}{u_2}
\arrow{sw,b}{f_2}\\ 
\node{}\node{H}
\end{diagram}
\] 
is commutative. Note that this diagram is obtained 
by reversing all arrows in the previous diagram 
and relabelling them. 

\begin{theorem} 
$((u_1, u_2), G_1\cup G_2))$ is the coproduct of 
$G_1$ and $G_2$ in $\textswab{S}$. 
\label{the:coproduct} 
\end{theorem} 

The \textit{cross product} of 
graphs $G_1$ and $G_2$ may be defined 
by the requirement that 
for every graph $G$ and homomorphisms 
$f_1 : G\rightarrow G_1$ and $f_2 : G\rightarrow G_2$, 
there exists a unique homomorphism 
$f : G\rightarrow G_1\times G_2$ 
so that the following diagram is commutative. 
\[
\begin{diagram} 
\node{}\node{G}\arrow{sw,t}{f_1}\arrow{s,l}{f}\arrow{se,t}{f_2}\\ 
\node{G_1}\node{G_1\times G_2}\arrow{w,t}{p_1}\arrow{e,t}{p_2}\node{G_2} 
\end{diagram}
\] 
where $p_1$ and $p_2$ are natural projections 
(homomorphisms). 

The \textit{disjunction} of 
graphs $G_1$ and $G_2$ may be defined 
by the requirement that 
for every graph $G$ and homomorphisms 
$f_1 : G\rightarrow G_1$ and $f_2 : G\rightarrow G_2$, 
there exists a unique homomorphism 
$f : G\rightarrow G_1\vee G_2$ 
so that the following diagram is commutative. 
\[
\begin{diagram} 
\node{}\node{G}\arrow{sw,t}{f_1}\arrow{s,l}{f}\arrow{se,t}{f_2}\\ 
\node{G_1}\node{G_1\vee G_2}\arrow{w,t}{p_1}\arrow{e,t}{p_2}\node{G_2} 
\end{diagram}
\] 
where $p_1$ and $p_2$ are natural projections. 

The \textit{strong product} of 
graphs $G_1$ and $G_2$ is defined to be 
the union of their cross and cartesian 
products. For definitions of products, 
the reader may also refer to a recent 
monograph \cite{imrichk:2000}. 

Many binary graph operations are interpreted 
categorically in \cite{kilpk:2001,knauer:2006}, 
where the following were among results established 
there. 
\begin{enumerate} 
\itemsep 0pt\parskip 0pt 
\item 
The cross product $G_1\times G_2$ with projections 
is a product of $G_1$ and $G_2$ in $\textswab{Gra}$. 
\item 
The strong product $G_1\boxtimes G_2$ with projections 
is a product of $G_1$ and $G_2$ in $\textswab{EGra}$. 
\item 
The disjunction with projections 
is a product of $G_1$ and $G_2$ in $\textswab{CGra}$. 
\end{enumerate} 

These capture the essence of the products 
concerned in the respective categories. 

Categories $\textswab{Gra}$, $\textswab{CGra}$ and 
$\textswab{Egra}$ also have coproducts and tensor 
products \cite{kilpk:2001,knauer:2006}. 
It was also shown in \cite{kilpk:2001,knauer:2006} 
that products and coproducts in these three 
categories have right adjoints. 

\section{Invariants} 
\label{sec:invar} 
Investigations about invariants in various fields 
of mathematics always concern the action of a group 
(usually a subgroup of the automorphism group). 
Combinatorics is not an exception. 
Let $\mathscr{G}$ be a set of graphs and $S$ be a set. 
A mapping $f : \mathscr{G}\rightarrow S$ is called an 
\textit{invariant} of graphs if for all $G, H\in \mathscr{G}$, 
$G\simeq H\Rightarrow f(G)\simeq f(H)$. Of course, 
for $S$ a set of numbers or a set of sequences of 
numbers, the second $\simeq$ is just $=$. In terms of mappings, 
a function taking its argument as a graph $G$ is an invariant if for 
each automorphism $\varphi$ of $G$, $f(\varphi (G)) = f(G)$, 
or simply, $f\varphi = f$, as the above is true for all graphs 
$G$ in the given family. Thus, a graph invariant 
may be presented by the following diagram. 
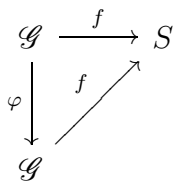
\[
\begin{diagram} 
\node{\mathscr{G}}\arrow{s,t}{\varphi}\arrow{e,t}{f}\node{S}\\ 
\node{\mathscr{G}}\arrow{ne,t}{f} 
\end{diagram}
\] 
Taking into account the condition 
$G\simeq H\Rightarrow f(G)\simeq f(H)$, 
a graph invariant for a family of graphs 
may be presented also by the diagram 
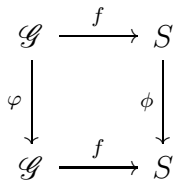
\[
\begin{diagram} 
\node{\mathscr{G}}\arrow{s,t}{\varphi}
\arrow{e,t}{f}\node{S}\arrow{s,t}{\phi}\\ 
\node{\mathscr{G}}\arrow{e,t}{f}\node{S} 
\end{diagram}
\] 
where $\varphi$ is a graph isomorphism and 
$\phi$ is an isomorphism of $S$. 

For example, (1) if $S$ 
is the set of all integer sequences, then the degree sequence 
function $f = {\bf d}$ is an invariant, since for each 
automorphism $\varphi$ of $G$, ${\bf d}\varphi = {\bf d}$; 
(2) if $S =\textswab{G}$ is the category of all groups then the automorphism 
group function $f =\mathrm{Aut}$ is an invariant since 
obviously $\mathrm{Aut}\varphi =\mathrm{Aut}$; 
(3) the determinant of the adjacency matrix is another 
example of an integer invariant; 
(4) the spectrum of a graph is an example of invariants; 
so also is the largest eigenvalue. 

A subgraph $H\subseteq G$ is called an \textit{invariant} 
subgraph if $\text{Aut}(G)(H) = H$. 

If $S\subseteq\mathbb{R}$ then the invariant 
$f: \mathscr{G}\rightarrow S$ is called a graph \textit{parameter}. 
In particular, any integer valued invariant is an 
example of a graph parameter. These include, of course, 
the order, size, diameter, girth, circumference, connectivity, 
edge connectivity, independence number, covering number, chromatic 
number, edge chromatic number and Ramsey number.  A significant 
part of graph thery dedicates itself to the study of graph 
parameters. 

\section{Algebraic Objects} 
\label{sec:mathobj} 
We have already stated that graphs, groups, rings, fields, and 
partially ordered sets and hence lattices are instances 
of binary and ternary relational systems. 
We shall consider more formally in this section 
some other algebraic systems. 

The specialization begins from the most abstract concept, 
the concept of a category.  For the formal axioms 
for categories, the reader may see \cite{maclane:1971}, 
where the following is explicitly stated and established. 

\begin{proposition} 
Every category is a graph. 
\label{prop:category} 
\end{proposition} 

As we have considered groups already, and category 
theory was essentially born out of a deep connection 
between groups and topological spaces, topological 
spaces will be addressed now. 

Every topological space is a hypergraph, and hence 
every topological space is a relational system. 

\begin{proposition} 
Every module is a ternary relational system. 
\label{prop:modules} 
\end{proposition} 

\textsl{Proof}: By definition, a module is an abelian 
group $M$ (a ternary relational system as seem above) 
together with a ring homomorphism 
$f : R\rightarrow {\rm End}(M, M)$. 
Now $f$ is a ternary relation over $M$. Hence a module 
is a ternary relational system. \pfend  

Thus, every vector space is a ternary relational 
system, and every algebra is a ternary relational system. 

From category theory, topology and algebra, we now 
return to combinatorics and consider matroids. 
Since every matroid is a hypergraph, 
every matroid is a relational system. 

Whereas the concept of a category captures 
mathematical concepts from algebro-axiomatic point 
of view (see \cite{maclane:1971} pages 10-12), 
relational systems capture them in 
an elementary combinatorial way. 
Having said about the generality achieved 
by the concept of relational systems, 
it needs to be pointed out that this 
generality is useful only as a proper 
generalization, as the concept of 
operations is considerably more special than 
that of a relation. Operations certainly 
bear more properties and these are 
exploited in the study of algebraic objects 
such as groups, rings, fields and modules. 

\section{Inductive Sets} 
\label{sec:recurs} 
A partially ordered set is said to be \textit{well founded} 
if every descending chain is finite (this is 
the Jordan-Dedekind descending chain condition). 
A well founded partial order is also abbreviated as a 
\textit{well founded order}. 
As a special type of partial order, 
every well founded order is a graph. 

An \textit{inductive class} $\mathscr{C}$ is usually 
understood as a set of objects 
such that a subset $\mathscr{B}\subseteq\mathscr{C}$ is 
designated and for each $X\in\mathscr{C}\backslash\mathscr{B}$ 
there is a well defined reduction $\rho$ such that 
$\rho (X)\in\mathscr{C}$. But the following 
definition is more essential. 

\begin{definition}{\rm 
A set $S$ is called \textit{inductive} if 
there is a well founded order on $S$. 
} 
\label{def:indset} 
\end{definition} 

The set of minimal elements is the set 
$\mathscr{B}$ in the previous paragraph. 
Note that nothing is said about whether the 
set of minimal elements is \textsl{finite}. 
In fact, consider the set of all positive 
integers excepting $1$, under the binary relation of 
divisibility: $a\leq b\Leftrightarrow a\;|\; b$. 
This relation is a well founded order, 
and the set of minimal elements (the set 
of all primes) is infinite. This is also 
an example of a well founded order that is 
not a well quasi order (as defined at 
the end of this section). 

The following statement says something about 
the domain of usage of the important principle 
of mathematical induction. 

\begin{proposition} 
The mathematical induction principle 
is valid on a set $S$ if and only if 
$S$ has a well founded order. 
\label{prop:mindprinciple} 
\end{proposition} 

This is one of the most basic statements in discrete mathematics. 
However, in classrooms this principle was 
taught in a way that gives an impression that 
this is a review of junior highschool mathematics. 
The importance of the fact that this principle 
should be understood here as a \textsl{characterization} 
or \textsl{complete determination} of the nature 
of a set on which mathematical induction 
may be used, is usually ignored or misconveyed! 
The true implication of this statement 
is that if a set $S$ has a well founded order 
then the mathematical induction may be applied, 
and if mathematical induction may be applied 
to elements of a set $S$ then $S$ has a well founded order. 
Unfortunately, this is seldom 
done. In an inductive set $S$, if every antichain 
is finite, then $S$ is called \textit{finitely generated}. 
Note also that each inductive set is a graph. 

\begin{definition}{\rm 
A reflexive and transitive binary relation is called 
a \textit{quasi-ordering}. A quasi-ordering $\leq$ on 
a set $S$ is a \textit{well quasi-ordering}, and 
the elements of $S$ are \textit{well quasi-ordered} 
by $\leq$, if for every infinite sequence 
$x_0, x_1, \cdots\in X$, there exist indices $i < j$ 
such that $x_i\leq x_j$. 
} 
\label{def:wqo} 
\end{definition} 

\begin{proposition} {\rm (\cite{diestel:1997} page 252)} 
A quasi-ordering is a well quasi-ordering 
if and only if every antichain is finite and 
every descending is finite. 
\label{prop:wqoiff} 
\end{proposition} 

Since a quasi order is a binary relation, 
every well quasi-ordering is an oriented graph. 
No one can deny the importance of well quasi orders 
in the theory of graphs. It is, however, usual 
to encounter a denial of the importance of 
a well founded order in the theory of graphs. 

Since a partially ordered set is 
a graph according to our definition, 
order preserving mappings between two 
partially ordered sets, and more specifically 
sets with well founded orders, is nothing but 
a graph homomorphism between the oriented graphs. 

\section{Contractions and Minors} 
\label{sec:contrmin} 
It was mentioned in Section \ref{sec:catgraphs} 
that a contraction is a preconnected egamorphism. 
An equivalent formulation is by way of a connected 
partition of $V(G)$. 

Let $G = (V, E)$ be a graph. 
For $X, Y\subseteq V(G)$, denote 
\[(X, Y) = \{xy : x\in X, y\in Y, xy\in E(G)\}.\]

A \textit{contraction} of $G$ is defined 
to be a partition $\{V_1, V_2, \cdots , V_s\}$ 
of $V$ such that for each $i = 1, 2, \cdots , s$, 
the induced subgraph $G|_{V_i}$ is connected. 
This partition gives rise to a natural mapping 
(this is a preconnected egamorphism) from $G$ to a graph $H$, 
also called a \textit{contraction} (graph) of $G$. 
The contraction (graph) $H$ is the graph with  
\[V(H) =\{V_1, V_2, \cdots , V_s\},\quad  
E(H) =\{V_iV_j : i\ne j,\; (V_i, V_j)\ne\emptyset\}.\] 
The mapping $f$ is called a \textit{contraction} (mapping) 
(or preconnected egamorphism) from $G$ 
onto $H$, and $G$ is said to be \textit{contractible} to $H$. 

The graph $K_1$ is a contraction of any connected graph $G$ 
since $\{V\}$ is a partition of $V$ and $G = G|_V$ is connected. 
Any automorphism of $G$ is a contraction 
since it is a permutation of the trivial partition 
of $V$ into single vertices. In particular, $1 : G\rightarrow G$ 
is a contraction. 

Suppose that $R\subseteq G$ is a connected subgraph. 
Then the contraction of $R$ in $G$, denoted $G/R$, 
is given by the partition 
\[\left\{V(R), \{v_1\}, \cdots , \{v_m\}\right\}\] 
where $V(G) - V(R) = \{v_i : 1\leq i\leq m\}$. 

Let $r\geq 1$ be an integer and denote 
by $e^r = (uv)^r$ the presence of $r$ parallel 
edges between vertices $u$ and $v$ in a multigraph. 
A contraction $f: G\rightarrow H$ is called 
a \textit{faithful} contraction if 
\[E(H) = \{(V_iV_j)^r : |(V_i, V_j)| = r, i\ne j\}.\] 
For an undirected graph, $(V_i, V_j)$ may be 
typed as $[V_i, V_j]$ or just $V_iV_j$. 

A contraction may be understood in various ways, but 
this shall not be our concern here in this paper. It is 
only noted here that this definition is adequate for 
directed graphs and infinite graphs, and is 
equivalent to stating that $f : G\rightarrow H$ 
is a preconnected egamorphism. 

A graph $H$ is a \textit{minor} of $G$, 
if $G$ has a subgraph contractible to $H$. 
That is, there is a subgraph $K\subseteq G$ 
and a contraction $f : K\rightarrow H$. 
This is the same as saying that the following 
diagram commutes. 
\[ 
\begin{diagram} 
\node{K}\arrow{s,b}{f}\arrow{r,t}{\eta}\node{G}\\ 
\node{H}\arrow{ne,b}{\mu} 
\end{diagram} 
\] 

This diagram may be used to prove some elementary 
properties of minor inclusions. First, since 
$1 : H\rightarrow H$ is a contraction, we have 
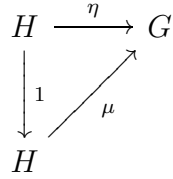
\[ 
\begin{diagram} 
\node{H}\arrow{s,b}{1}\arrow{r,t}{\eta}\node{G}\\ 
\node{H}\arrow{ne,b}{\mu} 
\end{diagram} 
\] 
Hence, 

(1) $H\subseteq G\Rightarrow H\leq G$. 

As a simple consequence, $G\leq G$ 
(the reflexivity of the binary relation $\leq$.) 

If $f : G\rightarrow H$ is a contraction, then 
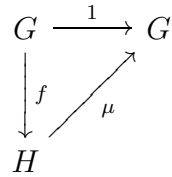
\[ 
\begin{diagram} 
\node{G}\arrow{s,b}{f}\arrow{r,t}{1}\node{G}\\ 
\node{H}\arrow{ne,b}{\mu} 
\end{diagram} 
\] 
Hence, we have proved 

(2) If $f: G\rightarrow H$ is a contraction, then $H\leq G$. 

Note that the converse is not true in general. 
For an example, $K_{3,3}\leq P$ where $P$ is 
the Petersen graph, but there is no contraction 
$f : P\rightarrow K_{3,3}$ (prove this!) 
 
Denote by $\dot{G}$ a subdivision (i.e., a homeomorph) 
of a graph $G$. As a corollary to (2), we have 

(3) If $\dot{H}\subseteq G$ then $H\leq G$.  

Now, the transitivity of the binary relation 
$\leq$ may also be established by using 
the diagram defining a minor. 

(4) $(J\leq H)\wedge (H\leq G)\Rightarrow J\leq G$; 

\textsl{Proof}: 
Consider the diagram 
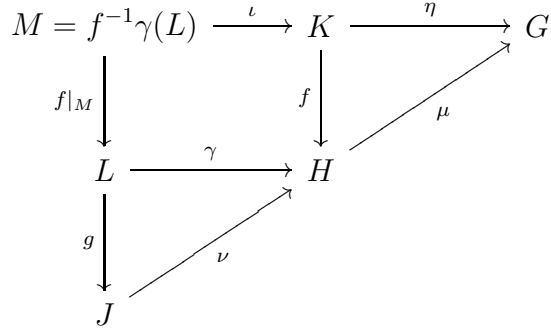
\[
\begin{diagram} 
\node{M=f^{-1}\gamma (L)}\arrow{e,t}{\iota}\arrow{s,l}{f|_{M}}\node{K}
\arrow{s,l}{f}\arrow{e,t}{\eta}\node{G}\\ 
\node{L}\arrow{e,t}{\gamma}\arrow{s,l}{g}\node{H}
\arrow{ne,r}{\mu}\\ 
\node{J}\arrow{ne,r}{\nu} 
\end{diagram}
\] 
In this diagram, $M = f^{-1}\gamma (L)\subseteq K\subseteq G$ 
is a subgraph of $G$, $f\iota =\gamma f|_M$, 
$f|_M(M) = f|_Mf^{-1}\gamma (L) =\gamma (L)\simeq L$ 
and $gf|_M : M\rightarrow J$ is a contraction since 
composition of contractions is a contraction. 
Hence $J\leq G$ by the definition of a minor. 
\pfend 

We cite, without proof, 
two further elementary properties of 
minor inclusions. 

(5) If $\Delta (G)\leq 3$, then 
$H\leq G\Leftrightarrow\dot{H}\subseteq G$; 

(6) If $H\leq G$ and $G$ is planar then $H$ is also planar. 

It will now be proved that the binary relation 
of minor inclusion is very close to being 
antisymmetric for a family of finite graphs. 
\begin{proposition} 
Let $G$ and $H$ be finite graphs. 
If $H\leq G$ and $G\leq H$ then $G\simeq H$. 
\label{prop:asymminor} 
\end{proposition} 

\textsl{Proof}: Suppose that $G$ and $H$ 
are finite graphs and that $H\leq G$ 
and $G\leq H$. Then we have diagrams 
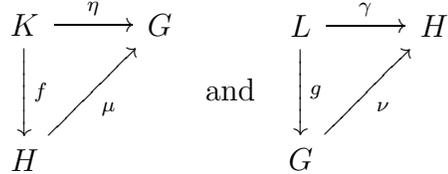
\[ 
\begin{diagram} 
\node{K}\arrow{s,b}{f}\arrow{r,t}{\eta}\node{G}\\ 
\node{H}\arrow{ne,b}{\mu} 
\end{diagram} 
\;\;\textrm{and}\;\; 
\begin{diagram} 
\node{L}\arrow{s,b}{g}\arrow{r,t}{\gamma}\node{H}\\ 
\node{G}\arrow{ne,b}{\nu} 
\end{diagram} 
\] 
Since $|G|\leq |L|\leq |H|$ and $|H|\leq |K|\leq |G|$, 
we have $|G| = |H|$. Hence $G$ is a spanning 
subgraphs of $H$ and $H$ is a spanning subgraph 
of $G$. (In particular, $V(G) = V(H)$.) 
This means. $E(H)\subseteq E(G)$ and $E(G)\subseteq E(H)$. 
Hence, $E(G) = E(H)$. Hence, 
\[f : G\rightarrow H,\; g: H\rightarrow G\] 
are both contractions and bijections
\[f : V(G)\rightarrow V(H),\; g: V(H)\rightarrow V(G).\] 
Since $f$ and $g$ are contractions, hence 
$uv\in E(G)$ if and only if $f(u)f(v)\in E(H)$ and 
hence $G\simeq H$. 
\pfend 

Thus, if isomorphic graphs are regarded as 
equal (which is a reasonable agreement), 
then the binary relation on finite graphs 
provided by minor inclusion is a partial order. 
Hence, a well quasi order in a family of 
graphs is essentially a well founded order. 
This, in addition to the principle of 
mathematical induction as stated in this note, 
is our reason why well founded order 
is interesting. 
 
\section{Transformation Graphs} 
\label{sec:transfg} 
In Definition \ref{def:graphasbinrel}, 
allow $V$ to be a specific set of mathematical objects. 
Then $E$ is a binary relation between mathematical objects, 
which may be given by a well defined set of transformations. 
Then Definition \ref{def:graphasbinrel} 
itself becomes the definition of a 
\textit{transformation graph}. The difference in the two definitions 
is that the set $V$ is an abstract set in the former 
and it is a specific set in the latter case. 
The binary relation $E$ is a binary relation 
defined on an abstract set in the former and 
it is one between specific mathematical objects 
in the latter case. 

Let us now consider special examples of transformation 
graphs that are sufficiently general and important 
in mathematics. 

(1) Consider the category of all finite groups 
and group homomorphisms. This is a transformation graph, 
which might have as well been called homomorphism graph 
of finite groups. The study of this graph comprises 
an essential part of the theory of groups. 

(2) Consider the category of all topological spaces 
and continuous mappings. This is a transformation graph. 
The study of this graph comprises an essential part of topology. 
The study of the homomorphism from the graph of (1) 
to the graph of (2) includes homotopy and homology. 

(3) Let $V$ be any finite set of 
positive integers, and let $E$ be the binary relation 
of divisibility: for $a, b\in V$, $(a, b)\in E$ if 
$a$ divides $b$. This tansformation graph has been 
called the \textit{divisibility graph} 
in \cite{chartrandmsz:2001}. 
 
(4) Let $G$ be a connected finite simple graph. 
Let $V$ be the set of its different 
spanning trees. For spanning trees $T_1, T_2\in V$, 
let $(T_1, T_2)\in E$ if $|\|T_1\| -\|T_2\|| = 2$. 
This is the \textit{tree transformation graph} 
of $G$ which had been studied by Whitney as early as 
1927 (I found no direct reference to this 
in the collecttion in my vicinity). 
The binary relation by which the edges 
are defined in the tree transformation graph was known 
in the literature as the \textsl{fundamental exchange} 
of edges. 
 
(5) Let $G$ be any finite simple graph with a perfect 
matching. Let $V$ be the set of its different 
perfect matchings. For matchings $M_1, M_2\in V$, 
let $(M_1, M_2)\in E$ if there is a unique 
$(M_1, M_2)$-alternating cycle $C$ in $G$. This 
may be called the \textit{matching transformation graph} 
of $G$. 
 
(6) Let $V$ be the set of all perfect matchings 
in a hexagonal system (see \cite{zhanggc:1988,zhangz:1995}), 
and let $E$ be the binary relation where $(M_1, M_2)\in E$ 
for $M_1, M_2\in V$ if a hexagon is $(M_1, M_2)$-alternating. 
This gives the concept of a $Z$-transformation 
graph \cite{zhangz:1995} of perfect matchings in a hexagonal system. 

(7)  Let $\mathbf{d}$ be a graphic degree sequence, 
and $\mathscr{R}(\mathbf{d})$ be the set 
of all isomorphism classes of finite simple graphs 
with degree sequence $\mathbf{d}$. 
For $G, H\in\mathscr{R}(\mathbf{d})$, 
$(G, H)\in E(\mathscr{R}_\mathbf{d})$ 
if there exist $ab, cd\in E(G)$ 
with $ac, bd\not\in E(G)$ such that 
\[H = (G - \{ab, cd\})\cup\{ac, bd\}.\] 
This is called the \textit{realization} graph 
of $\mathbf{d}$. This graph has been studied 
for many interesting parameters in \cite{punnim:2005}. 

(8) This is a new concept, some special cases of which 
have been studied recently. Let $r$ be a fixed 
positive integer, $H$ be a fixed graph 
and $\cdot$ be a fixed binary operation. Then for 
a graph $G$, an $(H, r)$-\textsl{transformation} graph 
$J = T_{H, r}(G)$ may be defined by assigning 
\[V(J) =\{S\subseteq E(G) : |S| = r\},\; 
E(J) =\{ST : S, T\in V(J), G|_{S\cdot T}\supseteq H\}.\] 
The transformation graph $T_{K_{1,2}, 1}(G) = L(G)$ is 
the usual \textsl{line graph}, and for 
$r = 2$ and $H = K_{1,2}$, the graphs $T_{H, r}(G)$ 
have been studied in \cite{lilz:2007}. 
The binary operation $\cdot$ being a natural 
product in an appropriate graph category 
seems to have not been of much attention. 

These examples point to the sources of transformation graphs: 
transformation graphs arise from 
(1) a set of mathematical objects; 
(2) a set of subobjects of a mathematical object. 

Tree transformation graphs have been studied in 
\cite{cummins:1966,mayedas:1965} 
where it was established that tree transformation graphs 
are connected. 
Line graphs and super line graphs have been studied in 
\cite{baggabv:1995,baggabv:1999,baggabv:2004,
beineke:1968,beineke:1971,chartrands:1969,zamfirescu:1970}. 
$Z$-transformation graphs of hexagonal systems 
have been studied extensively 
in \cite{cyving:1988,zhanggc:1988,zhangl:1995,zhangz:1995,zhangz:2000}. 
The study of hexagonal systems were directly motivated by 
organic chemistry. Matching transformation graphs have been 
studied in \cite{bau:2001,bauh:2003} where it was established 
that these graphs are $2$-connected. Euler tour graphs 
have been studied in \cite{li:1997,zhangg:1986}. 
Divisibility graphs 
have been studied in \cite{chartrandmsz:2001,gerasz:2003}. 
Transformation graphs based on some specific edge operations 
were studied in \cite{goddards:1996}. 
Switching transformation graphs 
or realization graphs have been investigated in 
\cite{eggletonh:1979,hakimi:1962,havel:1955,punnim:2005,taylor:1982}.  
Oriented transformation graphs of the quasi order 
arising from minor inclusion has been intensively 
investigated by Robertson and Seymour 
(see Diestel \cite{diestel:1997}, Chapter 12 for a sketch). 

Problems of connectivity of transformation graphs 
and those of traversals have been investigated. 
Measures of compactness and metric properties 
of the graphs such as diameter have also been of interest. 
In general, fundamental properties 
(combinatorial, geometric, topological or algebraic) 
of transformation graphs are of interest. 
A more detailed report on transformation 
graphs will be given in another paper by the author. 

\thebibliography{99}
\itemsep=0pt\parskip 0pt 
{\small 
\bibitem{alons:1992} 
N. Alon and J.H. Spencer, 
The Probabilistic Method, 
Wiley, New York 1992. 
\bibitem{bafnabf:1999} 
V. Bafna, P. Berman and T. Fujito, 
A 2-approximation algorithm for the undirected feedback vertex set problem, 
\textit{SIAM J. Discrete Math.}, \textbf{12}(1999), 289-297. 
\bibitem{baggabv:1995} 
K.S. Bagga, L.W. Beineke and B.N. Varma, 
Super line graphs, in Y. Alavi, A. Schwenk (Eds.), 
Graph Theory, Combinatorics, and Applications, 
Wiley-Interscience, New York 1995, 35-46. 
\bibitem{baggabv:1999} 
K.S. Bagga, L.W. Beineke and B.N. Varma, 
Independence and cycles in super line graphs, 
\textit{Australas. J. Combin.}, \textbf{19}(1999), 171-178. 
\bibitem{baggabv:2004} 
K.S. Bagga, L.W. Beineke and B.N. Varma, 
Old and new generalizations of line graphs, 
\textit{J. Math. Math. Sci.}, \textbf{29}(2004), 1509-1521. 
\bibitem{bau:2001} 
S. Bau, The connectivity of matching transformation graphs 
of cubic bipartite plane graphs, \textit{Ars Comb}., 60(2001), 161-169. 
\bibitem{bauh:2003} 
S. Bau and M.A. Henning, 
Matching transformation graphs of cubic bipartite plane graphs, 
\textit{Discrete Math}., 262(2003), 27-36. 
\bibitem{beineke:1968} 
L.W. Beineke, 
Derived graphs and digraphs, 
in \textit{Beitr\"{a}ge zur Graphentheorie}, 
H. Sachs, H. Voss and H. Walther (Eds.) Teubner, Leipzig 1968, 17-33. 
\bibitem{beineke:1971} 
L.W. Beineke, 
Derived graphs and derived complements, 
in \textit{Recent Trends in Graph Theory}, 
M. Capobianco, J.B. Frechen, M. Krolik (Eds.) 
Springer Verlag, New York 1971, 15-24. 
\bibitem{birkhoff:1967} 
G. Birkhoff, Lattice Theory, 
Amer. Math. Soc. Colloq. Publ. 25, 
Providence R.I. 1967. 
\bibitem{bollobas:1985} 
B. Bollob\'{a}s, 
Random Graphs, Academic Press, London 1985. 
\bibitem{chartrands:1969} 
G. Chartrand and M.J. Stwewart, 
The connectivity of line graphs, 
\textit{Math. Ann.}, \textbf{182}(1969), 170-174. 
\bibitem{chartrandmsz:2001} 
G. Chartrand, R. Muntean, V. Saepholphat and P. Zhang, 
Which graphs are divisor graphs?, 
\textit{Congr. Num.}, 151(2001), 189-200. 
\bibitem{cummins:1966} 
R.L. Cummins, Hamilton circuits in tree graphs, 
\textit{IEEE Trans., Circuit Theory}, 13(1966), 82-90. 
\bibitem{cyving:1988} 
S.J. Cyvin and I. Gutman, Kekul\'{e} Structures in Benzenoid Hydrocarbons, 
Lecture Notes in Chemistry 46, Springer Verlag, Berlin 1988. 
\bibitem{diestel:1997} 
R. Diestel, 
Graph Theory, Springer Verlag, New York 1997. 
\bibitem{eggletonh:1979} 
R.B. Eggleton and D.A. Holton, 
Graphic sequences, 
Combinatorial Mathematics VI (Proc. Sixth Aistral. Conf., University of 
New England, Armidale 1978), 
Lecture Notes in Mathematics 748(1979), 1-10. 
\bibitem{gerasz:2003} 
R. Gera, V. Saenpholphat and P. Zhang, 
Divisor graphs with triangles, 
\textit{Congr. Numer.}, 161(2003), 19-32. 
\bibitem{goddards:1996} 
W. Goddard and H.C. Swart, 
Distances between graphs under edge operations, 
\textit{Discrete Math}., 161(1996), 121-132. 
\bibitem{godsil:1993} 
C.D. Godsil, 
Algebraic Combinatorics, Chapman \& Hall, New York 1993. 
\bibitem{godsilr:2001} 
C.D. Godsil and G.F. Royle, 
Algebraic Graph Theory, 
Springer Verlag, New York 2001. 
\bibitem{grosst:2001} 
J.L. Gross and T.W. Tucker, 
Topological Graph Theory, 
Dover, New York(?) 2001. 
\bibitem{hakimi:1962} 
S. Hakimi, 
On the realizability of a set of integers as the degree 
of the vertices of a graph, 
\textit{SIAM J. Appl. Math.}, \textbf{10}(1962), 496-506. 
\bibitem{havel:1955} 
M. Havel, 
A remark on the existence of finite graphs, 
\textit{Gasopis Pest. Mat.}, \textbf{80}(1955), 
477-480. (in Hungarian) 
\bibitem{helln:2004} 
P. Hell and J. Ne\v{s}et\v{r}il, 
Graphs and Homomorphisms, 
Oxford University Press 2004. 
\bibitem{imrichk:2000} 
W. Imrich and S. Klav\v{z}ar, 
Product Graphs, 
John Wiley \& Sons, New York 2000. 
\bibitem{kilpk:2001} 
M. Kilp and U. Knauer, 
Graph operations and categorical constructions, 
\textit{Acta. et Comment. Uni. Tartuensis de Math.}, 
5(2001), 43-57. 
\bibitem{kilpkm:2000} 
M. Kilp, U. Knauer and A.V. Mikhalev, 
Monoids, Acts and Categories, 
Walter de Gruyter, Berlin 2000. 
\bibitem{knauer:2006} 
U. Knauer, 
Semigroups and graph categories, 
(lecture notes from a series of lectures 
presented at Chiang Mai University, Thailand), 
2006. 
\bibitem{li:1997} 
X-L. Li, A lower bound for the connectivity of directed Euler tour 
transformation graphs, \textit{Discrete Math.}, 
\textbf{163}(1997), 101-108. 
\bibitem{lil:1999} 
D-M. Li and Y-P. Liu, 
A polynomial algorithm for finding the minimum feedback vertex set of a
3-regular simple graph, 
\textit{Acta Math. Sci.} (English Ed.), \textbf{19}(1999), 375-381. 
\bibitem{lilz:2007} 
H. Li, X. Li and H. Zhang, 
Path-comprehensive and vertex-pancyclic properties 
of super line graph $L_2(G)$, 
submitted. 
\bibitem{maclane:1971} 
S. Mac Lane, 
Categories for the Working Mathematician, Springer-Verlag, New York 1971. 
\bibitem{maurer:1973} 
S.B. Maurer, Matroid basis graphs, II, J. Comb. Theory, 14B(1973), 216-240. 
\bibitem{mayedas:1965} 
W. Mayeda and S. Seshu, 
Generation of trees without duplications, 
\textit{IEEE Trans. Circuit Theory}, 
12(1965), 181-185. 
\bibitem{mohart:2001} 
B. Mohar and C. Thomassen, 
Graphs on Surfaces, 
Johns Hopkins University Press 2001. 
\bibitem{punnim:2005} 
N. Punnim, 
Switchings, realizations, and interpolation theorems 
for graph parameters, 
\textit{Intern. J. Math. Math. Sci.}, 
\textbf{13}(2005), 2095-2117. 
\bibitem{taov:2006} 
T. Tao and V.H. Vu, 
Additive Combinatorics, 
Cambridge University Press, Cambridge 2006. 
\bibitem{taylor:1982} 
R. Taylor, 
Switchings constrained to 2-connectivity in simple graphs, 
\textit{SIAM J. Alg. Disc. Meth.}, 3(1982), 114-121.  
\bibitem{zamfirescu:1970} 
T. Zamfirescu, 
On the line connectivity of line graphs, 
\textit{Math. Ann.}, \textbf{187}(1970), 305-309. 
\bibitem{zhangg:1986} 
F-J. Zhang and X-F. Guo, Hamilton cycles in Euler tour graphs, 
\textit{J. Comb. Theory.}, \textbf{1B}(1986), 1-8. 
\bibitem{zhanggc:1988} 
F-J. Zhang, X-F. Guo and R-S. Chen, Z-transformation graphs of perfect 
matchings of hexagonal systems, \textit{Discrete Math.}, 
\textbf{72}(1988), 405-415. 
\bibitem{zhangl:1995} 
F-J. Zhang and X-L. Li, Hexagonal systems with forcing edges, 
\textit{Discrete Math.}, \textbf{140}(1995), 253-263. 
\bibitem{zhangz:1995} 
F-J. Zhang and H-P. Zhang, A new enumeration method for Kekul\'{e} structures 
of hexagonal systems with forcing edges, 
\textit{J. Molecular Structure (Theochem)},
\textbf{331}(1995), 255-260. 
\bibitem{zhangz:2000} 
H-P. Zhang and F-J. Zhang, Plane elementary bipartite graphs, 
\textit{Discrete Appl. Math.}, \textbf{105}(2000), 291-311. 
} 
\end{document}